\documentclass[a4paper,12pt]{amsart}
\usepackage{amssymb}

\def\Fbar{\bar{F}}
\def\Gbar{\bar{G}}

\def\hbar{\bar{h}}

\def\Rarr#1{\buildrel #1\over \longrightarrow}

\def\iso{\buildrel \sim\over\to}

\def\Gb{{\mathfrak{b}}}
\def\Gg{{\mathfrak{g}}}
\def\Gh{{\mathfrak{h}}}
\def\Gm{{\mathfrak{m}}}

\def\Db{{\mathcal D}^b}
\def\CA{{\mathcal{A}}}
\def\CB{{\mathcal{B}}}
\def\CC{{\mathcal{C}}}
\def\CD{{\mathcal{D}}}
\def\CE{{\mathcal{E}}}

\def\CH{{\mathcal{H}}}

\def\CK{{\mathcal{K}}}
\def\CL{{\mathcal{L}}}

\def\CO{{\mathcal{O}}}
\def\CP{{\mathcal{P}}}

\def\BC{{\mathbf{C}}}

\def\BF{{\mathbf{F}}}
\def\BG{{\mathbf{G}}}
\def\BH{{\mathbf{H}}}

\def\BL{{\mathbf{L}}}

\def\BP{{\mathbf{P}}}
\def\BQ{{\mathbf{Q}}}

\def\BS{{\mathbf{S}}}

\def\BW{{\mathbf{W}}}

\def\BZ{{\mathbf{Z}}}

\def\Bs{{\mathbf{s}}}
\def\Bt{{\mathbf{t}}}

\def\Bw{{\mathbf{w}}}

\def\eps{\varepsilon}

\def\acyc{\operatorname{acyc}\nolimits}

\def\can{{\mathrm{can}}}

\def\mcomod{\operatorname{\!-comod}\nolimits}
\def\Comp{\operatorname{Comp}\nolimits}

\def\mdgMod{\operatorname{\!-dgmod}\nolimits}
\def\mdgperf{\operatorname{\!-dgperf}\nolimits}

\def\en{{\mathrm{en}}}
\def\End{\operatorname{End}\nolimits}

\def\Frob{\operatorname{Frob}\nolimits}
\def\Fun{\operatorname{Fun}\nolimits}

\def\Hom{\operatorname{Hom}\nolimits}

\def\id{\operatorname{id}\nolimits}
\def\im{\operatorname{im}\nolimits}

\def\Lie{\operatorname{Lie}\nolimits}
\def\mMod{\operatorname{\!-mod}\nolimits}

\def\mModgr{\operatorname{\!-modgr}\nolimits}

\def\opp{{\operatorname{opp}\nolimits}}

\def\mproj{\operatorname{\!-proj}\nolimits}

\def\Res{\operatorname{Res}\nolimits}

\def\Specm{\operatorname{Specm}\nolimits}

\def\ie{{\em i.e.}}

\def\tF{{\tilde{F}}}
\def\tG{{\tilde{G}}}

\def\tT{{\tilde{T}}}

\def\teps{{\tilde{\eps}}}
\def\teta{{\tilde{\eta}}}

\def\tPhi{{\tilde{\Phi}}}

\def\tUpsilon{{\tilde{\Upsilon}}}

\def\CHom{{{\mathcal H}om}}

\newtheorem{thm}{Theorem}[section]
\newtheorem{lemma}[thm]{Lemma}

\newtheorem{prop}[thm]{Proposition}

\newtheorem{conj}[thm]{Conjecture}

\theoremstyle{definition}
\newtheorem{rem}[thm]{Remark}

\input xypic
\xyoption{all}
\begin{document}
\author{Rapha\"el Rouquier}
\address{Rapha\"el Rouquier\\
Institut de Math\'ematiques de Jussieu --- CNRS\\
UFR de Math\'ematiques, Universit\'e Denis Diderot\\
2, place Jussieu\\
75005 Paris, FRANCE}
\email{rouquier@math.jussieu.fr}
\title{Categorification of the braid groups}
\date{Septembre 2004}
\maketitle

\section{Introduction}
Actions of braid groups on triangulated categories are quite widespread. They
arise for instance in representation theory, for constructible sheaves on
flag varieties and for coherent sheaves on Calabi-Yau varieties
(cf \cite{RouZi} and \cite{SeiTho} for early occurrences). In this work,
we suggest that not only the self-equivalences are important, but that
the morphisms between them possess some interesting structure.

\smallskip
Let $W$ be a Coxeter group and $\CC$ a triangulated category. We
consider gradually stronger actions of $W$ or its braid group $B_W$~:
\begin{enumerate}
\item[(i)] $W$ acting on $K_0(\CC)$
\item[(ii)] a morphism from $B_W$ to the group of isomorphism classes
of invertible functors of $\CC$
\item[(iii)] an action of $B_W$ on $\CC$.
\end{enumerate}

We construct a strict monoidal category $\CB_W$ categorifying
(conjecturally) $B_W$ and we propose an even stronger form of action~:

\begin{enumerate}
\item[(iv)] a morphism of monoidal categories
$\CB_W\to \CHom(\CC,\CC)$
\end{enumerate}

\smallskip
The first section is devoted to a construction of a self-equivalence
of a triangulated category, generalizing various constructions in
representation theory and algebraic geometry. This should be viewed as
a categorification of an action of $\BZ/2$.

In section \S\ref{sec2braid}, we construct a monoidal category
categorifying (a quotient of) the braid group.
It is a full subcategory of a homotopy
category of complexes of bimodules over a polynomial algebra. The setting
here is that of Soergel's bimodules.

Section \S\ref{BGG} is devoted to the category $\CO$ of a semi-simple
complex Lie algebra. There are classical functors that induce an action of
type (i).
We show how to use the constructions of \S\ref{sec2braid}, via results of
Soergel, to get a genuine action of the braid group and even the stronger
type (iv).

The case of flag varieties is considered in \S\ref{flag}. There again,
there is a classical action up to isomorphism of the braid group
on the derived category of constructible sheaves (type (ii)).
Using a result of Deligne and checking some compatibilities
for general kernel transforms (Appendix in \S\ref{compatkernel}),
one gets a genuine action of the braid
group (type (iii)). Now, using the link with modules over the
cohomology ring, we get another proof of this and even the stronger (iv).

\smallskip
In a work in preparation we study presentations by generators and
relations, homological vanishings and relation with the cohomology of
Deligne-Lusztig varieties.

\medskip
A few talks have been given in 1998--2000 on the main results of this work
(Freiburg, Paris, Yale, Luminy) and I apologize for the delay in putting them
on paper.

I would like to thank I.~Frenkel, M.~Kashiwara, M.~Khovanov and W.~Soergel
for useful discussions.

\section{Self-equivalences}
We describe a categorification of the notion of reflection with respect
to a subspace. The ambient space is $K_0$ of a triangulated category,
the subspace is another triangulated category and the embeddings and
projections are given by functors. We actually allow an automorphism
of the ``subspace'' category, which allows to categorify the $q$-analog
of a reflection.

\smallskip
We present here how a functor from a given triangulated category gives
rise to a self-equivalence of that category, when the functor satisfies
some conditions. Then, we give three special ``classical'' cases.
The first one concerns
constructible sheaves on a $\BP^1$-fibration (it occurs typically
with flag varieties, cf \S\ref{flag}). The second one 
deals with the case where the target category is the derived category of
vector spaces, where we recover the theory of spherical objects
and twist functors (it arises as counterparts of Dehn twists via mirror
symmetry). The last application essentially concerns
derived categories of finite dimensional algebras (it occurs in particular
within rational representation theory, cf \S\ref{BGG}).

\medskip
All functors between additive (resp. triangulated)
categories are assumed to be additive (resp. triangulated).

Given an additive category $\CC$, we denote by $K(\CC)$ the homotopy
category of complexes of objects of $\CC$.

Given an algebra $A$ over a field $k$,
we denote by $A\mMod$ the category of finitely
generated left $A$-modules. We put $A^\en=A\otimes_k A^\opp$, where
$A^\opp$ is the opposite algebra.

Given a graded algebra $A$, 
we denote by $A\mModgr$ the category of finitely generated graded $A$-modules. 

\subsection{A general construction}
\subsubsection{}
This section can probably be skipped in a first reading.

We will be working here with {\em algebraic} triangulated categories
(following Keller), a simple setting that provides functorial cones.

\smallskip
Let $\CE$ be a Frobenius category (an exact category with enough projective
and injective objects and where injective and projective objects coincide).
Let $\Comp_{\acyc}(\CE\mproj)$
be the category of acyclic complexes of projective objects of $\CE$.
Let $\Frob$ be the $2$-category of Frobenius categories,
with $1$-arrows the exact functors that send projectives to projectives and
$2$-arrows the natural transformations of functors.

The construction $\CE\mapsto \Comp_{\acyc}(\CE\mproj)$
is an endo-$2$-functor of $\Frob$.
The $2$-functor from $\CE$ to the $2$-category of triangulated
categories that sends $\CE$ to its stable category $\bar{\CE}$
factors through the previous functor.

The important point is that the category $\Comp_{\acyc}(\CE\mproj)$ has
functorial cones. Given $F,G:\Comp_{\acyc}(\CE\mproj)\to
\Comp_{\acyc}(\CE'\mproj)$ and $\phi:F\to G$,
then we have a well defined cone $C(\phi)$ of $\phi$ and
we have morphisms $G\to C(\phi)$ and $C(\phi)\to F[1]$ such that
$F\to G\to C(\phi)\to F[1]$ gives a distinguished triangle of functors
from $\bar{\CE}$ to $\bar{\CE'}$.

Note that if $\phi_0:F_0\to G_0$ is a morphism of functors (exact,
preserving projectives) between $\CE$ and $\CE'$, then we get via
$\Comp_{\acyc}(-)$ a morphism of functors
$\phi:F\to G$, with $F,G:\Comp_{\acyc}(\CE\mproj)\to
\Comp_{\acyc}(\CE'\mproj)$.

\smallskip
The category of functors (exact, preserving projectives)
$\Comp_{\acyc}(\CE\mproj)\to \Comp_{\acyc}(\CE'\mproj)$ is a Frobenius
category. We define the category $\mathrm{AlgTr}(\CE,\CE')$
to be its stable category. 
Its objects are the exact functors
$\Comp_{\acyc}(\CE\mproj)\to \Comp_{\acyc}(\CE'\mproj)$
and $\Hom_{\mathrm{AlgTr}(\CE,\CE')}(F,G)$ is the image
of $\Hom(F,G)$ in
$\Hom_{\Fun(\bar{\CE},\bar{\CE'})}(\Fbar,\Gbar)$.

This defines
the $2$-category of algebraic triangulated categories $\mathrm{AlgTr}$,
with objects the Frobenius categories $\CE$. We have a $2$-functor from
$\mathrm{AlgTr}$ to the $2$-category of triangulated categories obtained
by sending $\CE$ to $\bar{\CE}$. It is $2$-fully faithful.

\subsubsection{}
\label{constrequiv}
Let $\CC$ and $\CD$  be two algebraic triangulated categories,
$F:\CC\to\CD$, $G:\CD\to\CC$ be two functors and
$\Phi$ be a self-equivalence of $\CC$. Let there be given also
two adjoint pairs
$(F,G)$ and $(G,F\Phi)$. We have four morphisms (units and counits of
the adjunctions)
$$\eta:1_\CD\to F\Phi G,\ \ \ \eps:GF\Phi\to 1_\CC$$
$$\eta':1_\CC\to GF,\ \ \ \eps':FG\to 1_\CD.$$
Let $\Upsilon$ be the cocone of $\eps'$ and $\Upsilon'$ be the cone of
$\eta$~: there are distinguished triangles of functors
$\Upsilon\to FG\Rarr{\eps'} 1_\CD\rightsquigarrow$ 
and $1_\CD \Rarr{\eta} F\Phi G \to\Upsilon'\rightsquigarrow$.

Assume 
\begin{equation}
\label{GF}
1_\CC\xrightarrow{\eta'}GF\xrightarrow{\eps \Phi^{-1}}\Phi^{-1}
{\buildrel 0\over\rightsquigarrow}
\end{equation}
is a distinguished triangle.

\begin{prop}
\label{inversibilitetriang}
The functors $\Upsilon$ and $\Upsilon'$ are inverse self-equivalences of $\CD$.
\end{prop}

\begin{proof}
Let $\gamma$ be the map $F\Phi G\to \Upsilon'$ in the triangle above,
\ie, we have the distinguished triangle
$1_\CD \Rarr{\eta} F\Phi G\xrightarrow{\gamma}\Upsilon'
\rightsquigarrow$.
We have a commutative diagram with horizontal and vertical distinguished
triangles
$$\xymatrix{
&& \\
&& FG\ar@{~>}[u] \\
FG \ar[rr]^-{FG\eta}\ar[urr]^-{\id} && FGF\Phi G \ar[rr]^-{FG\gamma}\ar[u]_-{F\eps G} &&
 FG\Upsilon' \ar@{~>}[r] & \\
&& F\Phi G\ar[u]^-{F\eta'\Phi G}
}$$
The octahedral axiom shows that
$(FG\gamma)\circ(F\eta'\Phi G):F\Phi G\iso FG\Upsilon'$
is an isomorphism.

We have a commutative diagram
$$\xymatrix{
F\Phi G\ar[rr]^-{F\eta'\Phi G}\ar[drr]_-{\id} && FGF\Phi G\ar[rr]^-{FG\gamma}
 \ar[d]^-{\eps'F\Phi G} && FG\Upsilon' \ar[d]^-{\eps'\Upsilon'}\\
&& F\Phi G \ar[rr]_-{\gamma} && \Upsilon'
}$$
The distinguished triangle $\Upsilon \Upsilon'\to FG\Upsilon'\xrightarrow{\eps'\Upsilon'} \Upsilon'
\rightsquigarrow$ gives a distinguished triangle
$\Upsilon \Upsilon'\to F\Phi G\xrightarrow{\gamma} \Upsilon'\rightsquigarrow$,
hence $\Upsilon \Upsilon'\simeq 1_\CD$.

The case of $\Upsilon'\Upsilon$ is similar --- note that the triangle (\ref{GF})
shows that $GF\simeq \id_\CD\oplus \Phi^{-1}$, hence
$\Phi$ commutes with $GF$.
\end{proof}

\begin{rem}
One sees easily that $\eta F+F\Phi\eta':F\oplus F\Phi\iso F\Phi GF$ and
$G \eta +\eta'\Phi G:G\oplus \Phi G\iso GF\Phi G$ are isomorphisms.
One can show that the requirement that
$\eta FG+F\Phi\eta'G:FG\oplus F\Phi G\to F\Phi GFG$ and
$FG \eta +F\eta'\Phi G:FG\oplus F\Phi G\to FGF\Phi G$ are isomorphisms, instead of
the stronger requirement that (\ref{GF}) is a distinguished triangle, is enough
to get Proposition \ref{inversibilitetriang}.
\end{rem}

\subsubsection{}
Let us recall a version of Barr-Beck's Theorem
(\cite[\S VI.7, exercice 7]{Mac}, \cite[\S 4.1]{De2}).

Let $\CC$ be a category. A comonad is the data of a functor
$H:\CC\to\CC$, of $c:H\to H^2$ and $\eps:H\to\id$ such that
$(cH)\circ c=(Hc)\circ c$ and $(\eps H)\circ c=(H\eps)\circ c$. Note
that $\eps$ is determined by $c$.

A coaction of $(H,c,\eps)$ on an object $M$ of $\CC$ is the data of
$\rho:M\to H(M)$ such that
$\eps(M)\circ\rho=\id_M$ and $c\circ\rho=F(\rho)\circ\rho$.
The category $(H,c,\eps)\mcomod$ has objects the pairs $(M,\rho)$ and
a morphism $(M,\rho)\to (M',\rho')$ is a morphism $f:M\to M'$
such that $\rho'f=H(f)\rho$.

\smallskip
Let $\CA$ and $\CB$ be two abelian (resp. algebraic triangulated categories),
$T:\CA\to\CB$ an exact functor (resp. a triangulated functor).
Assume $T$ has a right adjoint $U$. Put $H=TU$, denote by
$\eps:H\to \id_\CB$ and $\eta:\id_\CA\to UT$ the counit and unit
of adjunctions and let $c=T\eta U:H\to H^2$.

\smallskip
We have a functor $\tT:\CA\to (H,c,\eps)\mcomod$ given by
$M\mapsto (TM,T\eta(M))$.

\smallskip
The following Theorem is an easy application of Barr-Beck's general result
to abelian and triangulated categories.

\begin{thm}
If $T$ is faithful, then $\tT:\CA\iso (H,c,\eps)\mcomod$ is an equivalence.
\end{thm}

\medskip
We deduce from this Theorem that the category $\CC$, together with the
functors $F$, $G$ and the adjunctions, is determined by $\CD$, $\Theta=FG$ and
$c=F\eta'G:\Theta\to \Theta^2$. We view this as a categorical version of
the ``fixed points'' construction.

\begin{rem}
Let $V=K_0(\CD)$, $U=K_0(\CC)$, $f=[F]:U\to V$ and $g=[G]:V\to U$. Assume
$[\Phi]=\id_U$. Then,
$gf=2 \id_U$ and $\theta=[\Theta]:x\mapsto x-fg(x):V\to V$ is an
involution. One recovers $U$ (up to unique isomorphism) from
$\theta$ acting on $V$ as $V^{\theta}$.
\end{rem}

\subsection{Applications}
\subsubsection{}
\label{P1}
We consider schemes of finite type type over an algebraic closure
of a finite field $\BF_q$ (the case of complex algebraic varieties
is similar).
Let $\pi:X\to Y$ be a smooth projective morphism already defined over
$\BF_q$. Assume the geometric fibers are projective lines.

Let $\Lambda$ be a field of coefficients (=an extension of $\BQ_\ell$
for $\ell\!\not|\ q$ a prime number).
Put $\CC=D^b(Y)$ and $\CD=D^b(X)$ (bounded derived categories
of constructible sheaves of $\Lambda$-vector spaces). 
Take $F=\pi^*$, $G=R\pi_*$ and $\Phi=?(1)[2]$.
We have a canonical isomorphism
(projection formula) $?\otimes R\pi_*\Lambda_X\iso R\pi_*\pi^*$.
Via this isomorphism, $\eta'$ becomes $\id\otimes\eta'(\Lambda_Y)$
and $\eps$ becomes $\id\otimes t$, where $t:R\pi_*\Lambda_X\to
\Lambda_Y(-1)[-2]$ is the trace map (an isomorphism on $\CH^2$).

So, the triangle (\ref{GF}) is obtained from
the triangle $\Lambda_Y\xrightarrow{\eta'(\Lambda_Y)}R\pi_*\Lambda_X
\xrightarrow{t}\Lambda_Y(-1)[-2]\rightsquigarrow$ by applying $?\otimes$.
This is indeed a distinguished triangle, for it is so at geometric fibers.

Let $\CL$ be a relative ample sheaf for $\pi$ and $c\in H^2(X,\Lambda(1))$ be
its first Chern class. The hard Lefschetz Theorem states that the composition
$\Lambda_Y\xrightarrow{\eta'(\Lambda_Y)}R\pi_*\Lambda_X\xrightarrow{c}
R\pi_*\Lambda_X(1)[2]\xrightarrow{t(1)[2]}\Lambda_Y$ is an isomorphism.
It follows that the connecting map in the triangle above is zero.

Thus, we are in the setting of \S \ref{constrequiv} and
we get a self-equivalence of $D^b(X)$.

\smallskip
This can be also constructed as a kernel transform.
Let $\alpha,\beta:X\times_YX\to X$ be the first and second projections.
Let $i:\Delta X\to X\times_Y X$ be the closed immersion of the diagonal
and $j:Z\to X\times_Y X$ be the open immersion of the complement of
$\Delta X$. Denote by $\teps:1_{D^b(X\times_Y X)}\to i_*i^*$ and
$\teta:Rj_!j^*\to 1_{D^b(X\times_Y X)}$ the adjunction morphisms.
One checks easily that there is a commutative diagram where the rows
are distinguished triangles
$$\xymatrix{
\Upsilon\ar[r] & \pi^*R\pi_*\ar[r]^-{\eps'}\ar[d]^\sim & 
 1_{\CD}\ar[d]^\sim \ar@{~>}[r] & \\
R\beta_*Rj_!j^*\alpha^*\ar[r]_-{R\beta_*\teta\alpha^*} &
 R\beta_*\alpha^*\ar[r]_-{R\beta_*\teps\alpha^*}&
 R\beta_*i_*i^*\alpha^*\ar@{~>}[r] & 
}$$
where the middle vertical map is the base change isomorphism.

Denote by $p,q:Z\to X$ the first and second projections.
Then, $\Upsilon\simeq Rp_!q^*$ and $\Upsilon'\simeq Rp_*q^!$.

\subsubsection{}
Assume we are in the setting of \S \ref{constrequiv} with
$\CC=D^b(k\mMod)$ where $k$ is a field and the categories and functors
involved are $k$-linear. There is an integer $n$ such that
$\Phi=?[n]$. Let $E=F(k)$. Then, $F\simeq E\otimes_k ?$ and
$G\simeq R\Hom(E,?)$. The morphism $\eps$ comes from
$t:\Hom(E,E[n])\to k$.

The morphism $\eps$ is the counit of an adjoint pair
$(G,F\Phi)$ if and only if 
$\dim_k\bigoplus_i\Hom(E,M[i])<\infty$ for all $M\in\CD$
and 
$\Hom(E,M)\times\Hom(M,E[n])\to k,\ (f,g)\mapsto t(gf)$ is a perfect
pairing for all $M\in\CD$.

The triangle (\ref{GF}) is distinguished if
and only if
$0\to k\cdot \id \to \bigoplus_i \Hom(E,E[i]) \xrightarrow{t} k\to 0$
is an exact sequence.

In other words, $E$ is an $n$-spherical object and $\Upsilon$, $\Upsilon'$ are the
corresponding twist functors of Seidel and Thomas \cite[\S 2b]{SeiTho}.
So, the framework above corresponds exactly to the twist functor theory
when $\CC\simeq D^b(k\mMod)$.

\begin{rem}
The case $\CC=D^b(k^d\mMod)$ also leads to interesting examples.
\end{rem}

\begin{rem}
It would be interesting to see if the construction of \S\ref{constrequiv}
can be used to construct
automorphisms of derived categories of Calabi-Yau varieties corresponding,
via Kontsevich's homological mirror symmetry conjecture, to
graded symplectic automorphisms on the mirror associated
to Lagrangian submanifolds more complicated than spheres.
\end{rem}

\subsubsection{}
\label{invRi}
Let us consider here two abelian categories $\CA$ and $\CB$
and $\tF:\CA\to\CB$, $\tG:\CB\to\CA$ and $\tPhi$ a self-equivalence of
$\CA$.
We assume we have two adjoint pairs $(\tF,\tG)$ and $(\tG,\tF\tPhi)$. So,
we have four morphisms (units and counits of the two adjunctions)
$$\teta:1_\CB\to \tF\tPhi\tG,\ \ \ \teps:\tG\tF\tPhi\to 1_\CA$$
$$\teta':1_\CA\to \tG\tF,\ \ \ \teps':\tF\tG\to 1_\CB.$$

Let $\tUpsilon$ be the complex $0\to \tF\tG\Rarr{\teps'} 1_\CB\to 0$ and
$\tUpsilon'$ the complex $0\to 1_\CB \Rarr{\teta} \tF\tPhi\tG\to 0$ (with
$\tF\tG$ and $\tF\tPhi\tG$ in degree $0$).
We put $\CC=K(\CA)$ and $\CD=K(\CB)$ and we denote by $F$, $G$, etc...
the extensions of $\tF$, $\tG$, etc... to $\CC$ and $\CD$.

Assume $\CB$ is artinian and noetherian (every object is a finite extension
of simple objects). If we have the equality 
$[\tG\tF]=[\id]+[\tPhi^{-1}]$ as endomorphisms of $K_0(\CA)$
(or more generally, if
$[\tF\tPhi \tG\tF\tG]=[\tF\tG\tF\tPhi \tG]=[\tF\tG]+[\tF\tPhi\tG]$
in $\End(K_0(\CB))$),
then, the conclusion of Proposition \ref{inversibilitetriang} remains valid. 

Let us justify this, following ideas of Rickard \cite[\S 3]{Ri}.
There is an adjoint pair $(\Upsilon',\Upsilon)$, hence there is a map
$u:\id\to \Upsilon\Upsilon'$ that doesn't vanish on a non-zero object of $\CB$.
One shows  that
$\Upsilon\Upsilon'$ is homotopy equivalent to a complex of functors with only
one non-zero term, $R$, in degree $0$ and $R$ is an exact functor.
The assumption on classes shows that $[R]=[\id]$.
So, $R$ sends a simple object to itself, for a simple
object is characterized amongst objects of $\CB$ by its class in $K_0(\CB)$.
In particular, $u:\id\to R$ is an isomorphism on simple objects.
So, $u$ is an isomorphism.

\section{The $2$-braid group}
\label{sec2braid}
\subsection{Coxeter group action}
\subsubsection{}
Let $(W,S)$ be a Coxeter system (with $S$ finite) and
$V=\bigoplus_{s\in S}k e_s$ be the reflection representation of $W$
over a field $k$. We assume the representation is faithful (this is always
the case if the characteristic is $0$).
Given $s,t\in S$, we denote by $m_{st}$ the order of $st$.
We assume that $2m_{st}$ is invertible in $k$, for all
$s,t\in S$ such that $m_{st}$ is finite.
We denote by $\{\alpha_s\}_{s\in S}$ the dual basis of $\{e_s\}_{s\in S}$
(so that $\ker(s-\id)=\ker\alpha_s$ for $s\in S$).
Let $B_W$ be the braid group of $W$. This is the group generated by
$\BS=\{\Bs\}_{s\in S}$ with relations
$$\underbrace{\Bs\Bt\Bs\cdots}_{m_{st}\text{ terms}}\simeq
\underbrace{\Bt\Bs\Bt\cdots}_{m_{st}\text{ terms}}$$
for any $s,t\in S$ such that $m_{st}<\infty$

Let $A=k[V]$ be the algebra of polynomial functions on $V$. All
$A$-modules considered in this section are graded.

We will sometimes identify an object $M$ of $K^b(A^\en\mModgr)$ with
the corresponding endofunctor $M\otimes_A -$
of $K^b(A\mModgr)$. In particular, we will sometimes omit the
symbols $\otimes_A$ when taking tensor products of bimodules
for the sake of clarity.

\subsubsection{}
The action of $W$ on $V$ induces an action on $A$, hence on $A\mModgr$ and
on $D^b(A\mModgr)$~: the element $w\in W$ acts by
$A_w\otimes_A -$ where $A_w$ is the $(A,A)$-bimodule equal to $A$ as a left
$A$-module, with right action of $a\in A$ given by right multiplication by
$w(a)$. We have an isomorphism of $(A,A)$-bimodules,
$\id\otimes 1:A_w\iso k[\Delta_w]$, where $\Delta_w=\{(w(v),v)\}_{v\in V}
\subset V\times V$.

We have a canonical isomorphism $A_w\otimes_A A_{w'}\iso A_{ww'}$
given by multiplication. Let $x=(x_1,\ldots,x_m)$ and
$y=(y_1,\ldots,y_n)$ be sequences of elements of $S$ such that
$x_1\cdots x_m=y_1\cdots y_n=w$.
We denote
by $c_{x,y}:A_{x_1} \cdots A_{x_m}\iso A_{y_1}\cdots A_{y_n}$
the isomorphism obtained by composing the multiplication map
$A_{x_1} \cdots A_{x_m}\iso A_w$ with the inverse of the
multiplication map $A_{y_1} \cdots A_{y_n}\iso A_w$.

\subsection{Braid group action}
Let us now construct a non-obvious lift of the action of $W$ on $D^b(A\mModgr)$
to an action of $B_W$ on $K^b(A\mModgr)$.

\subsubsection{}
For $s\in S$, we define the complex of $(A,A)$-bimodules
$$F_s=F_{\Bs}=0\to A\otimes_{A^s}A\Rarr{\eps'_s} A\to 0$$
where $A$ is in degree $1$ and $\eps'_s$ is the multiplication.

Since $A=A^s\oplus A^s\alpha_s$, the morphism of $A^\en$-modules
$$A_s\to A\otimes_{A^s}A(1),\
  a\mapsto a\otimes\alpha_s- a\alpha_s\otimes 1$$
induces an isomorphism 
$$f_{\Bs}:A_s\iso F_s(1) \text{ in }D^b(A^\en\mModgr).$$

\subsubsection{}
For $w\in W$, let $\Delta_{\le w}=\bigcup_{w'\le w}\Delta_{w'}$ and
$D_w=k[\Delta_{\le w}]$. Note that $D_s=A\otimes_{A^s}A$ for $s\in S$.

Given $w'\le w$, we have a canonical quotient map $D_w\to D_{w'}$ given by
restriction of functions.
We have
$$\Hom(D_w,D_{w'})=
\begin{cases}
k\cdot\can & \text{ if } w'\le w \\
0 & \text{ otherwise}
\end{cases}
$$

\subsubsection{}
In the next lemma, $0\to L\to M\to 0$ denotes a complex with $L$ in degree $0$.

\begin{lemma}
\label{declength2}
Assume $W$ is a finite dihedral group, \ie, $\dim V=2$, $S=\{s,t\}$
and $m_{st}<\infty$. Let $x\in W$ such that $tx>x$.
Then,
\begin{itemize}
\item[(i)]
$D_s(0\to D_{tx}\xrightarrow{\can}D_x\to 0)\simeq
(0\to D_x(-1)\xrightarrow{\id} D_x(-1)\to 0)\oplus
(0\to D_{stx}\xrightarrow{\can} D_x\to 0)$.
\item[(ii)]
$F_sD_x\simeq
(0\to D_x\xrightarrow{\id} D_x\to 0)\oplus
(0\to D_x(-1)\to 0\to 0)$.
\end{itemize}
\end{lemma}

\begin{proof}
Let us recall some constructions and results of Soergel 
\cite[Lemma 4.5, Proposition 4.6 and their proofs]{Soe4}. 
Since $2m_{st}$ is invertible, then given $u,u'$ two distinct reflections
of $W$, we have $\ker(u+\id)\not=\ker(u'+\id)$.

Given $I$ an ideal of $A\otimes A$
invariant
under $s\times 1$, we put $((A\otimes A)/I)^+=((A\otimes A)/I)^{(s\times 1)}$.

Let $r$ be the reflection of $W$ such that
$rx<x$ and $rx\not<tx$. Then, $\Delta_{x}+\Delta_{rx}$ is a hyperplane
of $V\times V$ and let $\beta\in V^*\times V^*$ be a linear form with kernel
this hyperplane. Let $M$ (resp. $N$) be the $(A^s\otimes A)$-submodule of
$D_{tx}$ generated by the image of the elements $\beta$ (resp $1$) of
$A\otimes A$. Then, $D_{tx}=M\oplus N$, $M\simeq D_x^+(-1)$ and
$N\simeq D_{stx}^+$ as $(A^s\otimes A)$-modules.

\smallskip
Let $M'$ (resp. $N'$) be the $(A^s\otimes A)$-submodule of $D_x$
generated by $\alpha_s\otimes 1$ (resp. $1$).
Then, $D_x=M'\oplus N'$, $M'\simeq D_x^+(-1)$ and
$N'=D_x^+$ as $(A^s\otimes A)$-modules.
Denote by $p:D_x\to M'$ the projection.

\smallskip
Let us show now that
$\beta\not\in (V^*)^s\times V^*$.
Equivalently, we need to show that
$(\Delta_x+\Delta_{rx})\cap (k\alpha_s\times 0)=0$.
This amounts to proving that $\im(\id-r)\not=k\alpha_s$. But this holds, since
$r\not=s$.

\smallskip
Let us now come to our problem.
Since $\beta\not\in (V^*)^s\times V^*$, it follows that
the image of $\beta$ in $(A\otimes A)/(A^s\otimes A)$ is a generator
as $(A^s\otimes A)$-module.
Consequently, the restriction of
$pf:D_{tx}\to M'$ to $M$ is surjective, hence it is an
isomorphism (we denote by $f:D_{tx}\to D_x$ the canonical map).
$$\xymatrix{
M\ar@{^{(}->}[r] \ar@/^2pc/[rrr]^\sim& D_{tx}\ar@{->>}[r]^-f &
 D_x\ar@{->>}[r]^p & M' \\
&& A\otimes A\ar@{->>}[ul] \ar@{->>}[u] \ar@{->>}[r] &
 A\otimes A/A^s\otimes A \ar@{->>}[u] \\
&& A^s\beta A\ar@{->>}[uull] \ar@{->>}[ur]\ar@{^{(}->}[u]
}$$

\smallskip
Finally, the multiplication map $A\otimes_{A^s}D_y^+\iso D_y$
is an isomorphism for any $y\in W$ with $sy<y$. 
We have shown that the complex
$A\otimes_{A^s}(0\to D_{tx}\xrightarrow{\can}D_x\to 0)$
is isomorphic to the direct sum of the complex
$0\to D_x(-1)\xrightarrow{\id} D_x(-1)\to 0$ and
a complex
$D=0\to D_{stx}\xrightarrow{\phi} D_x\to 0$.
Note that $\phi=r\cdot\can$ for some $r\in k$ and we need to prove that
$r\not=0$. The complex $0\to D_{tx}\xrightarrow{\can}D_x\to 0$
has zero homology in degree $1$, hence the same is true for
$D$. It follows that $r\not=0$.

\medskip
Let us now prove the second assertion. The multiplication map
$A\otimes_{A^s}D_x^+\to D_x$ is an isomorphism.
Since $D_x=D_x^+\oplus M'$ and $M'\simeq D_x^+(-1)$, we obtain
the second part of the Lemma.
\end{proof}

\begin{prop}
\label{tressebimod}
Take $s\not=t\in S$ with $m_{st}<\infty$.
We have braid relations
$$\underbrace{F_sF_tF_s\cdots}_{m_{st}\text{ terms}}\simeq
\underbrace{F_tF_sF_t\cdots}_{m_{st}\text{ terms}}$$
in $K^b(A\otimes A)$.
\end{prop}

\begin{proof}
We have a decomposition $V=V_1\oplus V_2$ under
the action of $\langle s,t\rangle$, with 
$V_1=V^{\langle s,t\rangle}$.
For the $(A,A)$-bimodules involved in
the Proposition, the right and left actions of $k[V_1]$ are identical.
So, we get the Proposition for $V$ from the Proposition for $V_2$ by
applying the functor $k[V_1^*]\otimes_k -$. It follows we can assume
$\dim V=2$.
So, we assume $W$ is finite dihedral with $S=\{s,t\}$.
We put $s_+=s$ and $s_-=t$.

Let $m=m_{st}$ and consider $i\le m$ and $\eps\in\{+,-\}$.
Let $\sigma^\eps_i=s_\eps s_{-\eps}s_\eps\cdots$ ($i$ terms) and
$D^\eps_i=D_{\sigma^\eps_i}$. We put $D^\eps=D_{s_\eps}$.
Consider the simplicial scheme over $V\times V$~:
$$\Delta_1\rightrightarrows \Delta_{\le s_+}\coprod\Delta_{\le s_-}
 \rightrightarrows \Delta_{\le s_+s_-}\coprod\Delta_{\le s_-s_+}
 \rightrightarrows\cdots\rightrightarrows
 \Delta_{\le \sigma^\eps_{i-1}}\coprod\Delta_{\le\sigma^{-\eps}_{i-1}}
\to \Delta_{\le\sigma_i^\eps}$$
where the maps are the inclusions.

We now define $F_i^\eps$ as the complex of $(A,A)$-bimodules coming from
the structural complex of sheaves of this simplicial scheme~:
$$F_i^\eps=
0\to D_i^\eps\xrightarrow{\tiny
\left(\begin{matrix}+\\+\end{matrix}\right)}
 D_{i-1}^\eps\oplus D_{i-1}^{-\eps}\xrightarrow{\tiny
\left(\begin{matrix}+&-\\+&-\end{matrix}\right)}
D_{i-2}^\eps\oplus D_{i-2}^{-\eps}\to \cdots\to
D^+\oplus D^- \xrightarrow{\tiny
\left(\begin{matrix}+&-\end{matrix}\right)}
 D_1\to 0$$
where the sign denotes the multiple of the canonical map considered
(we put $D_i^\eps$ in degree $0$).
We have $H^r(F_i^\eps)=0$ for $r>0$, since
$\Delta_{\le \sigma_r^+}\cap \Delta_{\le \sigma_r^-}=
\Delta_{\le \sigma_{r-1}^+}\cup \Delta_{\le \sigma_{r-1}^-}$ and we have
an exact sequence
$$0\to k[\Delta_{\le \sigma_r^+}\cup \Delta_{\le \sigma_r^-}]
\xrightarrow{\tiny\left(\begin{matrix}+\\+\end{matrix}\right)}
k[\Delta_{\le \sigma_r^+}]\oplus k[\Delta_{\le \sigma_r^-}]
\xrightarrow{\tiny \left(\begin{matrix}+&-\end{matrix}\right)}
k[\Delta_{\le \sigma_r^+}\cap \Delta_{\le \sigma_r^-}]\to 0.$$

\smallskip
The complex $F_1^\eps$ is isomorphic to $F_{s_\eps}$.
We will now show by induction on $i$ that $F_{s_\eps}F_i^{-\eps}$
is homotopy equivalent to $F_{i+1}^\eps$ for $\eps=\pm$. This will prove
the Proposition, since $F_m^+\simeq F_m^-$.

\smallskip
Let us consider the complex $C=F_{s_\eps}F_i^{-\eps}$. This is the total
complex of the double complex
$$\xymatrix{
D^\eps D_i^{-\eps} \ar[r]\ar[d] & D^\eps D_{i-1}^{\eps}\oplus
 D^\eps D_{i-1}^{-\eps} \ar[r]\ar[d] & D^\eps D_{i-2}^{\eps}\oplus
 D^\eps D_{i-2}^{-\eps}\ar[r]\ar[d] & \cdots\ar[r] & D^\eps D_1\ar[d] \\
D_i^{-\eps} \ar[r] & D_{i-1}^{\eps}\oplus D_{i-1}^{-\eps}
\ar[r] & D_{i-2}^{\eps}\oplus D_{i-2}^{-\eps}\ar[r] & \cdots\ar[r] & D_1
}$$

By Lemma \ref{declength2}, the complex
$0\to D^\eps D_r^{-\eps}\xrightarrow{\can} D^\eps D_{r-1}^{\eps}\to 0$
is isomorphic to the direct sum of
$0\to D_{r-1}^\eps(-1)\xrightarrow{\id}D_{r-1}^\eps(-1)\to 0$ and of
$0\to D_{r+1}^\eps\xrightarrow{\can} D_{r-1}^\eps\to 0$. Also, the complex
$0\to D^\eps D_r^\eps\xrightarrow{\can}D_r^\eps\to 0$
is isomorphic to the direct sum of
$0\to D_r^\eps\xrightarrow{\id}D_r^\eps\to 0$ and of 
$0\to D_r^\eps(-1)\to 0\to 0$.
It follows that $C$ is homotopy equivalent to a complex
$$C'=0\to D_{i+1}^\eps \to D_i^{\eps}\oplus D_i^{-\eps}\to\cdots\to D_1\to 0,$$
where the maps remain to be determined.
Since $F_{s_\eps}$ has non zero homology only in degree $0$ and that homology
is free as a right $A$-module, it follows that the homology of $C$ vanishes
in degrees $>0$.

\smallskip
To conclude, we have to show that a complex $X$ with the same terms as 
$F_i^\eps$ and with zero homology in degrees $>0$ is actually isomorphic
to $F_i^\eps$.
We have
$$X=
0\to D_i^\eps\xrightarrow{\tiny
\left(\begin{matrix}a_i\\c_i\end{matrix}\right)}
 D_{i-1}^\eps\oplus D_{i-1}^{-\eps}\xrightarrow{\tiny
\left(\begin{matrix}a_{i-1}&b_{i-1}\\c_{i-1}&d_{i-1}\end{matrix}\right)}
D_{i-2}^\eps\oplus D_{i-2}^{-\eps}\to \cdots\to
D^+\oplus D^- \xrightarrow{\tiny
\left(\begin{matrix}c_1&d_1\end{matrix}\right)}
 D_1\to 0$$
where the coefficients are in $k$ and the maps are corresponding multiples of
the canonical maps.

Take $r\le i$ minimal such that there is an entry of
$\tiny \left(\begin{matrix}a_r&b_r\\c_r&d_r\end{matrix}\right)$ that vanishes.
Assume for example $c_r=0$. Then, $a_{r-1}a_r=0$, hence $a_r=0$.
We have $b_rc_{r+1}=d_rc_{r+1}=b_rd_{r+1}=d_rd_{r+1}=0$. If
$b_r=d_r=0$, then $X$ is the sum of the subcomplex with zero terms in
degrees $\le i-r$ and the subcomplex with zero terms in degrees $>i-r$.
Otherwise, $c_{r+1}=d_{r+1}=0$, hence $X$ splits as the direct
sum of the subcomplex
$\cdots\to D_{r+1}^\eps\oplus D_{r+1}^{-\eps}\to D_r^\eps\to 0$ 
and the subcomplex
$0\to D_r^{-\eps}\to D_{r-1}^\eps\oplus D_{r-1}^{-\eps}\to\cdots$.
Now, a morphism $D_r^{-\eps}\to D_{r-1}^\eps\oplus D_{r-1}^{-\eps}$
is never injective, for the support of the left term is strictly larger than
the support of the right term. Consequently, the complex $X$ has non-zero
homology in degree $i-r$, which is a contradiction. We have proven that
none of the coefficients $a_r,b_r,c_r,d_r$ can be zero.

Let $Z$ be the closed subvariety of the affine space of coefficients
$a_r,b_r,c_r,d_r$
that define a complex (\ie, 
$\tiny \left(\begin{matrix}a_r&b_r\\c_r&d_r\end{matrix}\right)
\left(\begin{matrix}a_{r+1}&b_{r+1}\\c_{r+1}&d_{r+1}\end{matrix}\right)=0$)
and let $Z^0$ be its open subset corresponding
to non-zero coefficients. We have an isomorphism $Z^0\iso (\BG_m)^{2i-1}$,
$h:(a_r,b_r,c_r,d_r)_r\mapsto (a_r,c_r)_r$.
The action of $(\BG_m)^{2i}$ on the terms of the complex induce an
action on $Z$. The corresponding action on $Z^0\simeq
(\BG_m)^{2i-1}$ has a unique
orbit. It follows that $X$ is isomorphic to $F_i^\eps$.
\end{proof}

\subsubsection{}
Let us define the complex of $A^\en$-modules
$$F_{\Bs^{-1}}=0\to A\Rarr{\eta_s} A\otimes_{A^s}A(1) \to 0$$
where $A$ is in degree $-1$ and 
$\eta_s(a)=a\alpha_s\otimes 1+a\otimes\alpha_s$.

\begin{lemma}
\label{inverse}
The complexes $F_{\Bs}$ and $F_{\Bs^{-1}}$ are inverse to each other in
$K^b(A^\en\mModgr)$.
\end{lemma}

\begin{proof}
Let $\CC=K^b(A^s\mModgr)$ and $\CD=K^b(A\mModgr)$. Let $F=A\otimes_{A^s} ?$,
$G=A\otimes_A?$ and $\Phi=?(1)$. The morphisms of $(A^s,A^s)$-bimodules
$$\eps_s:A(1)\to A^s, 1\mapsto 0\text{ and }\alpha_s\mapsto 1
\ \ \text{ and}\ \ \eta'_s:A^s\to A, 1\mapsto 1$$
together with $\eta_s$ and $\eps'_s$ previously defined give rise
to adjoint pairs $(F,G)$ and $(G,F\Phi)$.

We have a split exact sequence of $(A^s,A^s)$-bimodules
$$0\to A^s\xrightarrow{\eta'_s}A\xrightarrow{\eps_s}A^s\to 0,$$
hence
we deduce the Lemma from Proposition \ref{inversibilitetriang}.
\end{proof}

By Proposition \ref{tressebimod} and Lemma \ref{inverse}, 
we have already obtained an action ``up to isomorphism'' of $B_W$ on
$K^b(A)$~:

\begin{prop}
\label{actionuptoiso}
The map $s\mapsto F_{\Bs}$ extends to a 
morphism from $B_W$ to the group of isomorphism
classes of invertible objects of $K^b(A^\en\mModgr)$.
\end{prop}

\subsection{Rigidification}
The key point here is that the rigidification of the braid relations
at the homotopy category level is equivalent to the one at the derived
category level, where the problem is trivial, since we have a genuine
action of $W$.
\subsubsection{}
Consider the morphism of $A^\en$-modules $A\otimes_{A^s}A\to A_s$ that
sends $1\otimes 1$ to $1$. It induces a
quasi-isomorphism $F_{\Bs^{-1}}(-1)\iso A_s$. We denote its inverse
(a morphism in $D^b(A^\en\mModgr)$) by $f_{\Bs^{-1}}$.

\smallskip
Now, let $v\in B_W$ and $v=t_1\cdots t_m=u_1\cdots u_n$ be two
decompositions in elements of $\BS\cup\BS^{-1}$.
By Proposition \ref{actionuptoiso}, the invertible objects $F_{t_1}\cdots
F_{t_m}$ and $F_{u_1}\cdots F_{u_n}$ of $K^b(A^\en\mModgr)$
are isomorphic, hence $$\Hom_\Box(F_{t_1}\cdots F_{t_m},
F_{u_1}\cdots F_{u_n})\simeq \End_\Box(A)=k,$$
where $\Box\in\{K^b(A^\en\mModgr),D^b(A^\en\mModgr)\}$.
It follows that the canonical morphism
$$\Hom_{K^b(A^\en\mModgr)}(F_{t_1}\cdots F_{t_m},
F_{u_1}\cdots F_{u_n})\iso \Hom_{D^b(A^\en\mModgr)}(F_{t_1}\cdots F_{t_m},
F_{u_1}\cdots F_{u_n})$$
is an isomorphism.

So, we have
a unique isomorphism
$$\gamma_{t,u}\in \Hom_{K^b(A^\en\mModgr)}(F_{t_1}\cdots F_{t_m},
F_{u_1}\cdots F_{u_n})$$
such that the induced element in 
$\Hom_{D^b(A^\en\mModgr)}(F_{t_1}\cdots F_{t_m}, F_{u_1}\cdots F_{u_n})$
corresponds to
$$c_{(t_1,\ldots,t_m),(u_1,\ldots,u_n)}:
A_{t_1}\cdots A_{t_m}\iso A_{u_1}\cdots A_{u_n}$$
via the quasi-isomorphisms 
$f_{t_1}\cdots f_{t_m}$ and $f_{u_1}\cdots f_{u_n}$.

\smallskip
We now define 
$G_v$ as the limit of the functors $F_{t_1}\cdots F_{t_m}$, where
$t=(t_1,\ldots,t_m)$ runs over the decompositions of $v$ in 
$\BS\cup\BS^{-1}$, with the
transitive system of isomorphisms $\gamma_{t,u}$.

There are unique isomorphisms $m_{v,v'}:G_vG_{v'}\iso G_{vv'}$ for $v,v'\in B_W$
and $m_1:G_1\iso A$ in $K^b(A^\en\mModgr)$ that are compatible with
the isomorphisms $c_{t,u}$, in $D^b(A^\en\mModgr)$. So, we get the
following result~:

\begin{thm}
\label{actionbimod}
The family $(G_v,m_{v,v'},m_1)$ defines an action of $B_W$ on $K^b(A\mModgr)$.
\end{thm}

This means we have a monoidal functor from
\begin{itemize}
\item the strict monoidal
category with set of objects $B_W$, with only arrows the identity maps and
with tensor product given by multiplication
\item to the strict monoidal category of endofunctors of $K^b(A\mModgr)$.
\end{itemize}

\begin{rem}
\label{rightaction}
Using tensor products on the right, one obtains a right action of
$B_W$ on $K^b(A\mModgr)$. This action commutes trivially with the
left action of $B_W$, so, we have an action of $B_W\times B_W^\opp$
on $K^b(A\mModgr)$.
\end{rem}

\subsubsection{}
We denote by $\CB_W$ the full subcategory of $K^b(A^\en\mModgr)$
with objects the $G_v$ for $v\in B_W$. The product
$G_v\boxtimes G_{v'}=G_{vv'}$ provides $\CB_W$ with the
structure of a strict monoidal category.
Define $G_v^*$ as $G_{v^{-1}}$.

We have obtained our ``categorification'' of the braid group~:

\begin{thm}
The category $\CB_W$ is a strict rigid monoidal category. Its
``decategorification'' is a quotient of $B_W$.
\end{thm}

\begin{conj}
The decategorification of $\CB_W$ is equal to $B_W$.
\end{conj}

\begin{rem}
One can show that the conjecture is true in type $A_n$, as a consequence of
\cite[Corollary 1.2]{KhovSei}.
\end{rem}

\subsubsection{}
\label{coinvariant}
Let $C=A/(A\cdot A^W_+)$ be the coinvariant algebra. Then, we get by
restriction
of functors an action of $B_W$ on $K^b(C\mModgr)$ and on $K^b(C\mMod)$.
We get as well monoidal functors from $\CB_W$ to the category of
self-equivalences of $K^b(C\mModgr)$ or $K^b(C\mMod)$. Note that
we get also right actions, and this gives a monoidal functor from
$\CB_W\times \CB_W^\opp$ to the category of 
self-equivalences of $K^b(C\mModgr)$ or $K^b(C\mMod)$.

\begin{rem}
Let $\CC$ be the smallest
full subcategory of $(A\otimes A)\mModgr$ containing the objects
$A\otimes_{A^s}A$ and closed under finite direct sums, direct summands
and tensor products. This is a monoidal subcategory of $(A\otimes A)\mModgr$
which is a categorification of the Hecke algebra of $W$, according to
Soergel. The quotient $\bar{\CC}$
of $\CC$ by the smallest additive tensor ideal
subcategory containing the $A\otimes_{A^{\langle s,t\rangle}}A$, where
$s,t\in S$ and $m_{st}\not=\infty$, is a categorification of
the Temperley-Lieb quotient of the Hecke algebra.

When $W$ has type $A_n$, an action of $\bar{\CC}$ on
an algebraic triangulated category is the same as the data on an
$A_n$-configuration of spherical objects \cite[\S 2.c]{SeiTho}.
\end{rem}

\section{Principal block of a semi-simple complex Lie algebra}
\label{BGG}

\subsection{Review of category $\CO$}
\subsubsection{}
Let $\Gg=\Lie G$, $\Gh\subseteq \Gb$ a Cartan and a Borel subalgebra.
Let $\CO$ be the Bernstein-Gelfand-Gelfand category of finitely generated
$\Gg$-modules which are diagonalizable
for $\Gh$ and locally finite for $\Gb$.
Denote by $Z$ the center of the enveloping algebra $U$ of $\Gg$. 
Let $P\subset \Gh^*$ be the weight lattice, $Q\subset \Gh^*$ be the
root lattice, $R$ (resp. $R^+$) be the set of roots (resp. positive
roots) and $\Pi$ the set of simple roots.

\subsubsection{}
We have a decomposition
$\CO=\bigoplus_\theta\CO_\theta$, where $\CO_\theta$ is the
subcategory of modules with central character $\theta$.
Let $D$ be a duality on $\CO$ that fixes simple modules (up to isomorphism).

Let $\Delta(\chi)=U\otimes_{U(\Gb)}\BC_\chi$ be the Verma module associated to
$\chi\in\Gh^*$. It has a unique simple quotient $L(\chi)$. We denote a
projective cover of $L(\chi)$ by $P(\chi)$.
We put $\nabla(\chi)=D\Delta(\chi)$.

Consider the dot action of $W$ on $\Gh^*$,
$w\cdot \lambda=w(\lambda+\rho)-\rho$ (we denote by $\dot{W}$ the
group $W$ acting via the dot action on $\Gh^*$), where
$\rho$ is the half-sum of the positive roots.

Given $\lambda\in\Gh^*$, let $\xi(\lambda)$ be the character by which
$Z$ acts on $L(\lambda)$ and $\Gm_\lambda$ be its kernel, an element of
$\Specm Z$,
the maximal spectrum of $Z$. The morphism
$\Gh^*\to \Specm Z,\ \lambda\mapsto\Gm_\lambda$
induces an isomorphism $\Gh^*/\dot{W}\iso \Specm Z$, \ie, an
isomorphism of algebras $h:Z\iso A^{\dot{W}}$ where $A=\BC[\Gh^*]$.
The simple objects in $\CO_\theta$ are those $L(\lambda)$ with
$\xi(\lambda)=\theta$.

\subsubsection{}
Consider $B$ the set of intersections of orbits of
$\dot{W}$ and of $Q$ on $\Gh^*$.
For $d\in B$, we denote by $\CO_d$ (or by $\CO_\mu$ for
a $\mu\in d$) the thick subcategory of
$\CO$ generated by the $L(\lambda)$ for $\lambda\in d$.
Then, $\CO=\bigoplus_{d\in B}\CO_d$ is the decomposition of
$\CO$ into blocks.

\smallskip
Let $\Lambda\in \Gh^*/P$ and $\lambda\in\Lambda$.
We have a root system $R_\Lambda=\{\alpha\in R| 
\langle \lambda,\alpha^\vee\rangle \in \BZ\}$ with set of simple roots
$\Pi_\Lambda\subset R^+$, Weyl group
$W_\Lambda=\{w\in W | w(\lambda)-\lambda\in Q\}$ and set of simple
reflections $S_\Lambda$ (they depend only on $\Lambda$).
Note that $R_\Lambda=R$ if and only if $\Lambda=P$.
We define 
$$\Lambda^+=\{\lambda\in\Lambda | \langle \lambda+\rho,\alpha^\vee\rangle\ge 0
\text{ for all }\alpha\in \Pi_\lambda\}$$
$$\Lambda^{++}=\{\lambda\in\Lambda | \langle
\lambda+\rho,\alpha^\vee\rangle> 0
\text{ for all }\alpha\in \Pi_\lambda\}.$$
Then, $\Lambda^+$ is
a fundamental domain for the action of $\dot{W}_\Lambda$ on $\Lambda$.
The module $L(\lambda)$ is finite dimensional if and only if
$\lambda\in P^{++}$.

\subsubsection{}
We define a translation functor between $\CO_d$ and $\CO_{d'}$ when
$d,d'\in B$ are in the same $P$-orbit. Take $\Lambda\in \Gh^*/P$ and
$\lambda,\mu\in\Lambda^+$. Let $\nu$ be the only element in
$W(\mu-\lambda)\cap \Lambda^{++}$. Then, we define
$T_\lambda^\mu:\CO_\lambda\to\CO_\mu, M\mapsto pr_\mu(M\otimes L(\nu))$
where $pr_\mu:\CO\to\CO_\mu$ is the projection functor.
Since $-w_0\nu\in \Lambda^{++}$ and
$L(\nu)^*\simeq L(-w_0\nu)$, it follows that the functors
$T_\lambda^\mu$ and $T_\mu^\lambda$ are left and right adjoint to each
other.

Let $d\in B$ containing $0$. The corresponding block
$\CO_0=\CO_d$ is the principal block of $\CO$. Note that
$d=W\cdot 0$ is a regular $W$-orbit and we put $L(w)=L(w\cdot 0)$, etc...

For $s\in S$, we fix $\mu\in P^+$ with stabilizer $\{1,s\}$ in $\dot{W}$.
We put $T^s=T_0^\mu$ and $T_s=T_\mu^0$ and
$\Theta_s=T_sT^s:\CO_0\to\CO_0$.

\subsubsection{}
Let $F_\Bs=F_s$ be the complex of functors on $\CO_0$ given by
$0\to \Theta_s\Rarr{\eps'_s} \id\to 0$
where $\eps'_s$ is the counit of adjunction ($\id$ is in degree 
$1$).

Let $F_{\Bs^{-1}}=
0\to \id\xrightarrow{\eta_s} \Theta_s\to 0$, where $\eta_s$
is the unit
of the other adjunction. Then, Rickard \cite[Proposition 2.2]{Ri} proved that
$F_\Bs$ and $F_{\Bs^{-1}}$ are inverse self-equivalences of $K^b(\CO_0)$
(this follows from \S \ref{invRi} by the classical character calculation
$[T^sT_s]=2[\id]$).

It is easy and classical that the $F_s$ induce an action of $W$
on $K_0(\CO_0)$ (the reflection $s\in S$
acts as $[F_s]$). This realizes the regular representation of $W$. A
permutation basis for this action
is provided by $\{[\Delta(w)]\}_{w\in W}$.

It seems difficult to check directly that the $F_s$ satisfy the braid
relations.
Using the equivalence between $\CO_0$ and perverse sheaves on the flag variety,
this can be deduced from \S\ref{flag}.

\subsection{Link with bimodules}
\label{appli}
\subsubsection{}
We start by recalling results of Soergel \cite{Soe1,Soe2,Soe3}
relating the category $\CO$ to modules over the coinvariant algebra.

Let $\Lambda\in\Gh^*/P$. We denote by $C_\Lambda=A/(A\cdot A_+^{W_\Lambda})$
the coinvariant
algebra of $(W_\Lambda,S_\Lambda)$ and $p_\Lambda:A\to C_\Lambda$
the canonical surjection.
Let $\lambda\in \Lambda^+$.
We denote by $t_\lambda:A\to A$ the translation by $\lambda$,
given by $f\mapsto (z\mapsto f(z+\lambda))$
We have Soergel's Endomorphismensatz \cite[Endomorphismensatz 7]{Soe1}~:

\begin{thm}
The image of the composite morphism $Z\Rarr{h} A^{\dot{W}}\hookrightarrow
A\Rarr{t_\lambda} A\xrightarrow{p_\Lambda} C_\Lambda$
is $C^{W_\lambda}_\Lambda$ and the canonical morphism 
$Z\to \End(P(w_0\cdot\lambda))$ factors through this morphism
$Z\to C^{W_\lambda}_\Lambda$. The induced morphism
$\sigma_\lambda:C^{W_\lambda}_\Lambda\iso \End(P(w_0\cdot\lambda))$ is an
 isomorphism.
\end{thm}

Let us now recall Soergel's Struktursatz \cite[Struktursatz 9]{Soe1}~:

\begin{thm}
\label{struktursatz}
The functor $\Hom(P(w_0\cdot\lambda),-):\CO_\lambda\mproj\to C^{W_\lambda}_\Lambda\mMod$ is
fully faithful.
\end{thm}

Let $\mu\in\Lambda$ be regular (\ie, with trivial stabilizer in $W_\Lambda$).

There is an isomorphism $\phi:T_\lambda^\mu P(w_0\cdot\lambda)\iso
P(w_0\cdot\mu)$.
Any such isomorphism $\phi$
induces a commutative diagram \cite[Bemerkung p.431]{Soe1}
$$\xymatrix{
C^{W_\lambda}_\Lambda\ar[d]_{\sigma_\lambda}\ar[r]^{\text{inclusion}} &
 C_\Lambda\ar[d]^{\sigma_\mu}  \\
\End(P(w_0\cdot\lambda))\ar[r]^{\phi_*T_\lambda^0} & \End(P(w_0\cdot\mu))
}$$

This gives us an isomorphism, via the adjunction $(T_\lambda^\mu,T^\lambda_\mu)$~:
$$\Res_{C^{W_\lambda}_\Lambda}^{C_\Lambda} \Hom(P(w_0\cdot\mu),?)\iso
\Hom(T_\lambda^\mu P(w_0\cdot\lambda),?)\iso
\Hom(P(w_0\cdot\lambda),T^\lambda_\mu(?))$$
between functors $\CO_\mu\to C^{W_\lambda}_\Lambda\mMod$.
So, we have a commutative diagram, with fully faithful
horizontal functors
$$\xymatrix{
\CO_\mu\mproj\ar[d]_{T^\lambda_\mu} \ar[rrr]^-{\Hom(P(w_0\cdot\mu),?)} &&&
 C_\Lambda\mMod \ar[d]^{\Res} \\
\CO_\lambda\mproj \ar[rrr]^-{\Hom(P(w_0\cdot\lambda),?)} &&& C^{W_\lambda}_\Lambda\mMod
}$$

\subsubsection{}
From the last commutative diagram, we deduce
\begin{prop}
There is a commutative diagram with fully faithful horizontal arrows
$$\xymatrix{
K^b(\CO_0\mproj) \ar[rrr]^-{\Hom(P(w_0),-)} \ar[d]_{F_s} &&&
  K^b(C\mMod) \ar[d]^{F_s} \\
K^b(\CO_0\mproj) \ar[rrr]^-{\Hom(P(w_0),-)} &&& K^b(C\mMod)
}$$
\end{prop}

So, we deduce from Theorem \ref{actionbimod} the following~:
given $v\in B_W$ and $v=t_1\cdots t_m=u_1\cdots u_n$ two
decompositions in elements of $\BS\cup\BS^{-1}$, there
is an isomorphism
$F_{t_1}\cdots F_{t_m}\iso F_{u_1}\cdots F_{u_n}$ between functors on
$D^b(\CO_0)$ coming by restriction from the isomorphism between functors
on $K^b(A\mModgr)$. These form a transitive system of isomorphisms, \ie

\begin{thm}
The functors $F_s$ induce an action of $B_W$ on $D^b(\CO_0)$.
\end{thm}

More precisely,
\begin{thm}
There is a monoidal functor from $\CB_W$ to the category of
self-equivalences of $D^b(\CO_0)$ sending $G_s$ to $F_s$.
\end{thm}

\begin{rem}
One has a similar statement for the deformed category $\CO$.

Note that we deduce from \S \ref{coinvariant} that there is also a right
action of $B_W$ on $D^b(\CO_0)$. We leave it to the reader to check that
this corresponds to the actions using Zuckerman functors, or equivalently,
Arkhipov functors.

In the graded setting (mixed perverse sheaves for example),
the left and right actions of $B_W$ should  be swapped by the self-Koszul
duality equivalence, cf \cite{BerFreKhov} (and
\cite[Conjecture 5.18]{BeiGi} for an analog in the equivariant case).

Various constructions have been given of weak actions of braid groups
on $D^b(\CO_0)$, cf \cite{AnStr,Ar,KhomMaz,MazStr,Str}.
\end{rem}

\section{Flag varieties}
\label{flag}
\subsection{Classical results}
Let $G$ be a semi-simple complex algebraic group with Weyl group $W$.

Let $\BW=\{\Bw\}_{w\in W}$. The braid group $B_W$ of $W$ is isomorphic
to the group with set of generators $\BW$ and relations
$\Bw\Bw'=\Bw'' \text{ when } ww'=w'' \text{ and }l(w'')=l(w)+l(w')$.

\smallskip
Let $\CB$ be the flag variety of $G$.
We decompose
$$\CB\times\CB=\coprod_{w\in W}\CO(w)$$
into orbits for the diagonal $G$-action.
Consider the first and second projections
$$\xymatrix{
 & \CO(w) \ar[dl]_{p_w} \ar[dr]^{q_w} \\
\CB & & \CB
}$$
Then, we have a functor
$$F_{\Bw}=R(p_w)_! (q_w)^*:D^b(\CB)\to D^b(\CB)$$
where $D^b(\CB)$ is the derived category of bounded complexes of
constructible sheaves of $\BC$-vector spaces over $\CB$.

First and last projections induce an isomorphism
$$\CO(w)\times_{\CB}\CO(w')\iso \CO(ww') \text{ when } l(ww')=l(w)+l(w').$$
This induces an isomorphism (cf \S \ref{noyaux})
$$\gamma_{w,w'}:F_{\Bw} F_{\Bw'}\iso F_{\Bw\Bw'} \text{ when }
l(ww')=l(w)+l(w').$$

For $s\in S$, then $F_\Bs$ is obtained as in \S\ref{P1} for the
canonical morphism $\pi_s:\CB\to\CP_s$, where $\CP_s$ is the variety of
parabolic subgroups of type $s$. So, $F_\Bs$ is invertible, with
inverse $F_{\Bs^{-1}}=R(p_s)_*(q_s)^!$.
It follows that $F_{\Bw}$ is invertible for $w\in W$, with inverse
$F_{\Bw^{-1}}=R(p_w)_*(q_w)^!$, hence 
we get a morphism from $B_W$ to the group of isomorphism
classes of invertible functors on $D^b(\CB)$.

\subsection{Genuine braid group action}

We have a commutative diagram

$$\xymatrix{
F_xF_yF_z \ar[r]^{\gamma_{x,y}F_z} \ar[d]_{F_x\gamma_{y,z}} &
  F_{xy}F_z \ar[d]^{\gamma_{xy,z}} \\
F_x F_{yz} \ar[r]_{\gamma_{x,yz}} & F_{xyz}
}$$

for $x,y,z\in \BW$ such that $l(x)+l(y)+l(z)=l(xyz)$, by Theorem
\ref{cocycle}.

Let $b\in B_W$ and $b=t_1\cdots t_m=u_1\cdots u_n$ with
$u_i\in \BW\cup\BW^{-1}$. Applying braid relations and the 
corresponding isomorphisms $\gamma$, we get various isomorphisms
$F_{t_1}\cdots F_{t_m}\iso F_{u_1}\cdots F_{u_n}$.
By Deligne \cite{De3}, they are all equal. Let us denote by
$\gamma_{t_\cdot,u_\cdot}$ their common value.

We now define
$$\tF_b=\lim_{(t_1,\cdots,t_n)} F_{t_1}\cdots F_{t_n}$$
where $(t_1,\cdots,t_n)$ runs over the set of sequences of
elements of $\BW\cup\BW^{-1}$ such that $b=t_1\cdots t_n$ and where
we are using the transitive system of isomorphisms $\gamma_{t_\cdot,s_\cdot}$.

We have now the following result

\begin{thm}
The assignment $b\mapsto \tF_b$ defines an action of $B_W$ on $D^b(\CB)$.
\end{thm}

\begin{rem}
Deligne \cite{De3} defines a variety $\CO_b$ with 
two morphisms $p_b,q_b:\CO_b\to\CB$ for any $b\in B_W^+$.
Then, the action of
$b$ on $\Db(\CB)$ is given by ${p_b}_!q_b^*$.
\end{rem}

\subsection{Link with bimodules}

\subsubsection{}
Fix a Borel subgroup $B$ of $G$. We consider the setting of \S\ref{sec2braid}
with $k=\BC$ and $V^*$ the complexified character group of $B$. 
In this section, we will consider the algebra $A$ with double
grading, \ie, $V^*$ is in degree $2$.

\smallskip
Let $C\iso H^*(\CB,\BC)$ be the Borel isomorphism (send a character of $B$
to the Chern class of the corresponding line bundle) and denote by
$\beta$ its inverse.

Let $I$ be a subset of $S$, $W_I$ the subgroup of $W$ generated by $I$,
$W^I$ be the set of minimal right coset
representatives of $W/W_I$ and $P_I$ the parabolic subgroup of $G$
of type $I$ containing $B$. Put $\CP_I=G/P_I$. Denote by
$\pi_I:\CB\to\CP_I$ the canonical morphism. The 
map $\pi_I^*:\bigoplus_i\Hom(\BC_{\CP_I},\BC_{\CP_I}[i])\to
\bigoplus_i\Hom(\BC_{\CB},\BC_{\CB}[i])$ induces, via $\beta$,
an isomorphism $\beta_I:\bigoplus_i\Hom(\BC_{\CP_I},\BC_{\CP_I}[i])\iso
C^{W_I}$.

\subsubsection{}
Consider the full subcategory $D_\sigma^b(\CP_I)$ of $D^b(\CP_I)$ of complexes
whose cohomology sheaves are smooth along $B$-orbits.
Given $w\in W^I$, let $\CL_w$ be the perverse sheaf corresponding to
the intersection cohomology complex of $\overline{BwP_I/P_I}$.
Let $\CL_I=\bigoplus_{w\in W^I}\CL_w$.
The dg-algebra $R\End(\CL_I)$ is formal and let
$R_I=\bigoplus_i\Hom(\CL_I,\CL_I[i])$. We have an equivalence 
$\CL_I\otimes ?$ from the category $R_I\mdgperf$ of perfect
differential graded $R_I$-modules to $D_\sigma^b(\CP_I)$.

The functor $\bigoplus_i\Hom(\BC_{\CP_I},?[i]):D^b_\sigma(\CP_I)\to
C^{W_I}\mModgr$
restricts to a fully faithful functor on the full subcategory
containing the $\CL_I[i]$. So, we get a fully faithful functor
$R_I\mdgperf\to K(C^{W_I}\mdgMod)$, hence a fully faithful functor
$\BH_I:D^b_\sigma(\CP_I)\to K(C^{W_I}\mdgMod)$, where we denote by
$K(C^{W_I}\mdgMod)$ the homotopy category of differential graded
$C^{W_I}$-modules.

\smallskip
As in \S\ref{appli}, we get a commutative diagram
$$\xymatrix{
D^b_\sigma(\CB)\ar[r]^-{\BH}\ar[d]_{R\pi_{I*}} & K(C\mdgMod)\ar[d]^{\Res} \\
D^b_\sigma(\CP_I)\ar[r]_-{\BH_I} & K(C^{W_I}\mdgMod)
}$$
and we deduce

\begin{prop}
Let $s\in S$. There is a commutative diagram with fully faithful
horizontal arrows
$$\xymatrix{
D^b_\sigma(\CB)\ar[r]^-{\BH}\ar[d]_{F_s} & K(C\mdgMod)\ar[d]^{F_s} \\
D^b_\sigma(\CB)\ar[r]_-{\BH} & K(C\mdgMod)
}$$
\end{prop}

In particular, we get a monoidal functor from
$\CB_W$ to the category of self-equivalences of $D^b_\sigma(\CB)$.

\begin{rem}
We believe the monoidal functor above is the restriction of a functor with
values in $D^b(\CB)$.
\end{rem}

\section{Appendix~: associativity of kernel transforms}
\label{compatkernel}

\subsection{Classical isomorphisms}
\subsubsection{}
We consider here
\begin{itemize}
\item
schemes of finite type over a field of characteristic $p\ge 0$
and the derived category of constructible sheaves of $\Lambda$-modules,
where $\Lambda$ is a torsion ring with torsion prime to $p$ or
$\Lambda$ is a $\BQ_\ell$-algebra, for $l$ prime to $p$

or
\item
locally compact topological spaces of finite soft $c$-dimension and the
derived category of constructible sheaves of $\BC$-vector spaces.
\end{itemize}

We will quote results pertaining to either of the two settings above,
depending on the convenience of references.
The maps involved will be concatenations of canonical isomorphisms.

We denote a derived functor with the same notation as the original
functor~: we write
$\otimes$ for $\otimes^\BL$, $f_!$ for $Rf_!$, etc...

\subsubsection{}
Let $f:Y\to X$ and $g:Z\to Y$ be two morphisms.
There are canonical isomorphisms \cite[2.6.6 and 2.3.9]{KaScha}
$$(fg)_!\iso f_!g_!\ \text{ and } \ (fg)^*\iso g^*f^*.$$

These isomorphisms satisfy a cocycle property
(cf \cite[Th\'eor\`eme 5.1.8]{De1} for the case $(-)_!$)~:

\begin{lemma}
\label{comp!}
Consider $X_3\Rarr{w} X_2\Rarr{v} X_1\Rarr{u} X_0$. Then, the following
diagrams are commutative
$$\xymatrix{
w^*v^*u^* \ar[r]\ar[d] & w^*(uv)^*\ar[d]& &
   u_!v_!w_! \ar[r]\ar[d] & (uv)_! w_!\ar[d] \\
   (vw)^*u^* \ar[r] & (uvw)^*& &
   u_!(vw)_! \ar[r] & (uvw)_!
   }$$
\end{lemma}

We will take the liberty to identify the functors
$v^*u^*$ and $(uv)^*$ through the canonical isomorphism.

\vskip 0.3cm
There are canonical isomorphisms \cite[2.6.18]{KaScha}
$$f^*(-_1\otimes -_2)\iso (f^*-_1)\otimes (f^*-_2)\ \text{ and }\
(-_1\otimes -_2)\otimes -_3\iso -_1\otimes (-_2\otimes -_3)$$

We identify the bifunctors
$f^*(-_1\otimes-_2)$ and
$(f^*-_1)\otimes (f^*-_2)$ through the canonical isomorphism.
Given $A_i\in D^b(X)$, $i\in\{1,2,3\}$,
we identify
$(A_1\otimes A_2)\otimes A_3$ with
$A_1\otimes (A_2\otimes A_3)$ and we denote this object by
$A_1\otimes A_2\otimes A_3$.

\vskip 0.3cm

Let
$\xymatrix{
X'\ar[d]_{f'}\ar[r]^{g'} & X \ar[d]^f \\
S' \ar[r]_{g}            & S
}$
be a cartesian square. Then, there is the canonical base change
isomorphism \cite[2.6.20]{KaScha}~:
$$g^*f_!\iso f'_!{g'}^*.$$

We have a canonical isomorphism \cite[2.6.19]{KaScha}
$$-_1\otimes (f_! -_2)\iso f_!(f^*-_1\otimes -_2).$$

\subsection{Kernel transforms}
\label{noyaux}
\subsubsection{}
\smallskip
Let us define a $2$-category $\CK$. 
\begin{itemize}
\item The $0$-arrows are the varieties.
\item $1$-arrows~: $\Hom(X,Y)$ is the family of 
$(K,U)$ where $U$ is a variety over $Y\times X$ and $K\in D^b(U)$.
\item $2$-arrows~: $\Hom((K,U),(K',U'))$ is the set of
$(\phi,f)$ where $f:U\iso U'$ is an isomorphism of $(Y\times X)$-varieties
and $\phi:K\iso f^*K'$.
\end{itemize}

We define the composition of $1$-arrows.
Consider the following diagram where the square is cartesian
\begin{equation}
\label{diag3}
\xymatrix{
& & V\times_Y U \ar[dl]_{\beta}\ar[dr]^{\alpha}\\
& V \ar[dl]_{p_4}\ar[dr]^{p_3} & & U \ar[dl]_{p_2} \ar[dr]^{p_1}\\
Z & & Y & & X
}
\end{equation}
Let $K\in D^b(U)$ and $L\in D^b(V)$.
We put $L\boxtimes K=\beta^*L\otimes \alpha^*K$.
The composition $(L,V)(K,U)$ is defined to be
$(L\boxtimes K,V\times_Y U)$.

\smallskip
Let us consider now the diagram with all squares cartesian
$$\xymatrix{
& & & W\times_Z V\times_Y U \ar[dl]_{b}\ar[dr]^{a} \\
& & W\times_Z V \ar[dl]_{\delta}\ar[dr]^{\gamma} & &
    V\times_Y U \ar[dl]_{\beta}\ar[dr]^{\alpha}\\
& W \ar[dl]_{p_6}\ar[dr]^{p_5} & & V \ar[dl]_{p_4}\ar[dr]^{p_3} & &
   U \ar[dl]_{p_2} \ar[dr]^{p_1}\\
T & & Z & & Y & & X
}$$
and take $M\in D^b(W)$.
We have
$$(M\boxtimes L)\boxtimes K=b^*(\delta^*M\otimes\gamma^*L)\otimes
(\alpha a)^*K\iso (\delta b)^* M\otimes a^*(\beta^*L\otimes \alpha^*K)=
M\boxtimes (L\boxtimes K).$$
This provides the associativity isomorphisms for $\CK$. With our conventions,
we will write $M\boxtimes L\boxtimes K$ for the objects in the
isomorphism above.
It is straightforward to check that $\CK$ is indeed a $2$-category.

\subsubsection{}
We put
$\Phi_K=\Phi_K^{p_2,p_1}={p_2}_!(K\otimes p_1^* -): D^b(X)\to D^b(Y).$

Let $c_{L,K}:\Phi_L\Phi_K\iso \Phi_{L\boxtimes K}$
be defined as the composition

\begin{align*}
{p_4}_!(L\otimes p_3^* {p_2}_!(K\otimes p_1^* -)) &\to 
{p_4}_!(L\otimes \beta_!\alpha^*(K\otimes p_1^* -))\\
& \to
{p_4}_!\beta_!(\beta^*L\otimes \alpha^*(K\otimes p_1^* -)) \\
&\to (p_4\beta)_!((\beta^*L\otimes \alpha^*K)\otimes \alpha^*p_1^* -)\\
&\to
(p_4\beta)_!((\beta^*L\otimes \alpha^*K)\otimes (p_1\alpha)^* -).
\end{align*}

Let $(\phi,f)\in\Hom_\CK((K,U),(K',U'))$. We have a commutative diagram
$$\xymatrix{
& U\ar[ddl]_{p_2}\ar[ddr]^{p_1}\ar[d]_{f}^{\sim} \\
& U'\ar[dl]^{p'_2}\ar[dr]_{p'_1} \\
Y && X
}$$
and we define $\Phi(\phi,f)$ as the composition
$$p_{2!}(K\otimes p_1^* -)\iso
p_{2!}'f_!(f^*K'\otimes f^*p_1^{'*} -)\iso
p_{2!}'(K'\otimes p_1^{'*} -).$$

\begin{thm}
\label{cocycle}
$\Phi$ is a $2$-functor from $\CK$ to the $2$-category of triangulated
categories.

We have $c_{M\boxtimes L,K}\circ(c_{M,L}\Phi_K)=
c_{M,L\boxtimes K}\circ(\Phi_M c_{L,K})$,
 {\em i.e.}, 
the following diagram commutes~:
$$\xymatrix{
\Phi_M\Phi_L\Phi_K \ar[r] \ar[d] & \Phi_{M\boxtimes L}\Phi_K \ar[d] \\
\Phi_M \Phi_{L\boxtimes K} \ar[r] & \Phi_{M\boxtimes L\boxtimes K}
}$$
\end{thm}

\subsubsection{}
The next two Lemmas deal with composition of base change isomorphisms.

For the first Lemma, see \cite[Lemme 5.2.5]{De1}~:

\begin{lemma}
\label{compoh}
Let
$\xymatrix{
X_2 \ar[r]^{f_2}\ar[d]_{h_2} & X_1\ar[r]^{f_1}\ar[d]_{h_1} & X \ar[d]_h\\
S_2 \ar[r]_{g_2} & S_1\ar[r]_{g_1} & S
}$ be a diagram with all squares cartesian.
Then, the following diagram commutes
$$
\xymatrix{
(g_1g_2)^*h_!\ar[rr]\ar[d] & & {h_2}_! (f_1f_2)^* \\
g_2^*g_1^*h_! \ar[r] & g_2^* {h_1}_!f_1^* \ar[r]& {h_2}_!f_2^*f_1^*\ar[u]
}$$
\end{lemma}

The second Lemma is \cite[Lemme 5.2.4]{De1}~:

\begin{lemma}
\label{compov}
Let
$\xymatrix{
X'_2\ar[d]_{f'_2}\ar[r]^{g_2} & X_2 \ar[d]^{f_2} \\
X'_1\ar[d]_{f'_1}\ar[r]^{g_1} & X_1 \ar[d]^{f_1} \\
S' \ar[r]_{g}            & S
}$ be a diagram with all squares cartesian.
Let $A\in D^b(S')$.
Then, the following diagram commutes
$$
\xymatrix{
A\otimes g^* (f_1f_2)_!- \ar[d]\ar[r] & A\otimes (f'_1f'_2)_! g_2^* -
\ar[r] & (f'_1f'_2)_!((f'_1f'_2)^*A\otimes g_2^* -)\\
A\otimes g^* {f_1}_!{f_2}_!- \ar[d]&& 
 {f'_1}_!{f'_2}_!({f'_2}^*{f'_1}^*A\otimes g_2^*-) \ar[u]\\
A\otimes {f'_1}_!g_1^*{f_2}_!- \ar[r]&
{f'_1}_!({f'_1}^*A\otimes g_1^*{f_2}_!-) \ar[r]&
{f'_1}_!({f'_1}^*A\otimes {f'_2}_!g_2^*-) \ar[u]\\
}$$
\end{lemma}

\begin{lemma}
\label{corKunneth}
Let $f:Y\to X$ and $A,B\in D^b(X)$ and $C\in D^b(Y)$. Then, the following
diagram commutes
$$\xymatrix{
A\otimes B\otimes f_!C \ar[r]\ar[d] & f_!(f^*(A\otimes B)\otimes C)\ar[d] \\
A\otimes f_!(f^*B\otimes C)\ar[r] & f_!(f^*A\otimes f^*B\otimes C)
}$$
\end{lemma}

\begin{proof}
The corresponding statement for $f_!$ replaced by $f_*$ is easy, the key point
is that the composition $f^*\xrightarrow{f^*\eta}
f^*f_*f^*\xrightarrow{\eps f^*} f^*$ is the identity of $f^*$, where
$\eta$ and $\eps$ are the unit and counit of the adjoint pair
$(f^*,f_*)$. The Lemma follows easily from this (in the algebraic case,
we have only to check in addition the trivial case where $f$ is an open immersion thanks to the transitivity of Lemma \ref{compov}, whereas
in the topological case we use the embedding $f_!C\subset f_*C$
for $C$ injective).
\end{proof}

\begin{lemma}
\label{triple}
Let
$\xymatrix{
X'\ar[d]_{f'}\ar[r]^{g'} & X \ar[d]^f \\
S' \ar[r]_{g}            & S
}$
be a cartesian square. Let $A\in D^b(S)$ and $B\in D^b(X)$.
Then, the following diagram commutes
$$\xymatrix{
g^*A\otimes g^*f_!B \ar[r]\ar[d] & g^*A\otimes f'_!{g'}^*B \ar[r] &
 f'_! ({f'}^*g^*A\otimes {g'}^*B)\ar[r] &
 f'_!({g'}^*f^*A\otimes {g'}^*B)\ar[d]\\
g^*(A\otimes f_!B)\ar[rr] && g^*f_!(f^*A\otimes B)\ar[r] & 
 f'_!{g'}^*(f^*A\otimes B)
}$$
\end{lemma}

\begin{proof}
As in the previous Lemma, one reduces to proving the analog of the Lemma
with $?_!$ replaced by $?_*$. This follows then from the easily checked
commutativity of the two diagrams
$$\xymatrix{
{f'}^*g^*f_*\ar[r]\ar[d]& {f'}^*f'_*{g'}^*\ar[d] & & 
  g^*\ar[r]\ar[d] & g^*f_*f^*\ar[d] \\
{g'}^*f^*f_*\ar[r]& {g'}^* &&
  f'_*{f'}^* g^* \ar[r] & f'_*{g'}^* f^*
}$$
where we have used the units and counits of the adjoint pairs
$(f^*,f_*)$ and $({f'}^*,f'_*)$.
\end{proof}

\begin{proof}[Proof of the Theorem]

We will show the commutativity of the following diagram
$$\xymatrix{
& & \Phi_M^{p_6,p_5}\Phi_{L\boxtimes K}^{p_4\beta,p_1\alpha}\ar[d]_\zeta
  \ar[drr]^{c_{M,L\boxtimes K}}\\
\Phi_M\Phi_L\Phi_K\ar[urr]^{\Phi_M c_{L,K}}\ar[drr]_{c_{M,L}\Phi_K}
 & & \Phi_{\delta^*M}^{p_6\delta,\gamma} \Phi_{L\boxtimes K}^{\beta,p_1\alpha}
\ar[rr]^{c_{\delta^*M,L\boxtimes K}} & & \Phi_{M\boxtimes L\boxtimes K}\\
& & \Phi_{M\boxtimes L}^{p_6\delta,p_3\gamma}\Phi_K^{p_2,p_1}
  \ar[u]^\xi \ar[urr]_{c_{M\boxtimes L,K}}
}$$
where $\zeta$ is the composition
\begin{align*}
{p_6}_!\left(M\otimes p_5^*(p_4\beta)_!(\beta^*L\otimes \alpha^* K\otimes
   (p_1\alpha)^*-)\right)&\to
{p_6}_!\left(M\otimes p_5^*{p_4}_!\beta_!(\beta^*L\otimes \alpha^* K\otimes
   (p_1\alpha)^*-)\right)\\
&\to {p_6}_!\left(M\otimes \delta_!\gamma^*\beta_!(\beta^*L\otimes \alpha^* K
   \otimes (p_1\alpha)^*-)\right)\\
&\to {p_6}_!\delta_!\left(\delta^*M\otimes \gamma^*\beta_!(\beta^*L\otimes
   \alpha^* K \otimes (p_1\alpha)^*-)\right)\\
&\to (p_6\delta)_!\left(\delta^*M\otimes\gamma^*\beta_!(\beta^*L\otimes
   \alpha^* K\otimes (p_1\alpha)^*-)\right)
\end{align*}
and $\xi$ the composition

\begin{align*}
(p_6\delta)_!\left(\delta^*M\otimes\gamma^*L\otimes (p_3\gamma)^*{p_2}_!
  (K\otimes p_1^*-)\right)&\to
(p_6\delta)_!\left(\delta^*M\otimes\gamma^*L\otimes \gamma^*p_3^*{p_2}_!
  (K\otimes p_1^*-)\right)\\
&\to (p_6\delta)_!\left(\delta^*M\otimes\gamma^*(L\otimes p_3^*{p_2}_!
  (K\otimes p_1^*-))\right)\\
&\to (p_6\delta)_!\left(\delta^*M\otimes\gamma^*(L\otimes \beta_!\alpha^*
  (K\otimes p_1^*-))\right)\\
&\to (p_6\delta)_!\left(\delta^*M\otimes\gamma^*(L\otimes
   \beta_!(\alpha^* K\otimes \alpha^*p_1^*-))\right)\\
&\to (p_6\delta)_!\left(\delta^*M\otimes\gamma^*(L\otimes
   \beta_!(\alpha^* K\otimes (p_1\alpha)^*-))\right)\\
&\to(p_6\delta)_!\left(\delta^*M\otimes\gamma^*\beta_!(\beta^*L\otimes
   \alpha^* K\otimes (p_1\alpha)^*-)\right)
\end{align*}

Let $u$ and $v$ be the compositions
$$u:{p_6}_!(M\otimes p_5^*{p_4}_!-)\to
{p_6}_!(M\otimes \delta_!\gamma^*-)\to
{p_6}_!\delta_!(\delta^*M\otimes \gamma^*-)\to
(p_6\delta)_! (\delta^*M\otimes \gamma^*-)$$
and
$$v:L\otimes p_3^*{p_2}_! (K\otimes p_1^* -)\to
L\otimes \beta_!\alpha^* (K\otimes p_1^* -)\to
\beta_!(\beta^*L\otimes \alpha^* (K\otimes p_1^* -))\to
   \beta_! (\beta^*L\otimes \alpha^*K\otimes (p_1\alpha)^* -).$$

Then, one has trivially
$$\zeta(\Phi_M c_{L,K})=u(L\otimes p_3^*{p_2}_! (K\otimes p_1^* -))\circ
{p_6}_!(M\otimes p_5^*{p_4}_!v)=\xi(c_{M,L}\Phi_K).$$

\vskip 0.5cm
The equality $c_{M,L\boxtimes K}=c_{\delta^*M,L\boxtimes K}\zeta$
follows from Lemma \ref{compov} applied to $g=p_5$, $g_1=\gamma$,
$g_2=a$, $f_1=p_4$, $f_2=\beta$, $f'_1=\delta$, $f'_2=b$ and $A=M$
and from Lemma \ref{comp!} applied to $u=p_6$, $v=\delta$ and $w=b$.

\vskip 0.5cm
The equality $c_{M\boxtimes L,K}=c_{\delta^*M,L\boxtimes K}\xi$
follows from Lemma \ref{compoh} applied to
$f_1=\alpha$, $f_2=a$, $g_1=p_3$, $g_2=\gamma$, $h=p_2$, $h_1=\beta$
and $h_2=b$, from Lemma \ref{corKunneth} applied to 
$f=b$, $A=\delta^*M$, $B=\gamma^*L$ and $C=(\alpha a)^*(K\otimes p_1^*-)$
and from Lemma \ref{triple} applied to
$f=\beta$, $g=\gamma$, $f'=b$, $g'=a$,
$A=L$ and $B=\alpha^*(K\otimes p_1^*-)$.
\end{proof}

\end{document}